\documentclass[aps,amsmath,amssymb,reprint,groupedaddress,floatfix]{revtex4-2}

\usepackage[usenames]{color} 
\usepackage{amssymb, amsmath, amsfonts, amsthm} 
\usepackage[utf8]{inputenc} 
\usepackage[T1]{fontenc}
\usepackage{graphicx}  
\usepackage{float} 

\usepackage{mathrsfs}
\usepackage{verbatim}
\usepackage{enumitem}

\usepackage{hyperref}
\hypersetup{
    colorlinks=true,
    linkcolor=blue,
    filecolor=blue,      
    urlcolor=blue,
    citecolor=blue,
}
\newcommand*\p{\mathbb{P}}

\def\1{\mathbf 1}
\def\d{{\rm d}}
\newcommand{\e}{{\rm e}}

\def\R{\mathbb R}

\newcommand{\avg}[1]{\left \langle #1 \right \rangle}

\newcommand{\T}{\mathcal{T}}

\begin{document}

\title{Asymptotic behaviour for escaping an interval by L\'evy flights under stochastic resetting}

\author{Bartosz Zbik}
\email{bartosz.zbik@york.ac.uk}
\affiliation{Department of Mathematics, University of York, Heslington, York, YO105DD, UK}

\author{Zbigniew Palmowski}
\email{zbigniew.palmowski@pwr.edu.pl}
\affiliation{Department of Applied Mathematics,
Wroc{\l}aw University of Science and Technology, Wybrze\.ze Wyspia\'nskiego~27
50-370 Wroc{\l}aw, Poland}

\date{13th October 2026}

\begin{abstract}
The resetting of stochastic processes has become a key mechanism appearing in many research fields, such as searching for a target, the transmission control protocol, biology, chemistry, finance, insurance, etc., with the main interest concerning the first passage time (FPT) by the process under resetting.
In this paper, we focus on the exit time $\T_L$ from the interval $(-L, L)$ of L\'evy flights (LF) under resetting, but we also consider spectrally negative L\'evy processes (SNLP).
As the main result, we find the asymptotics of the mean first passage time (MFPT) as $L$ tends to $+ \infty$.
For LF with $\alpha \in (0, 2)$, the scaling behaves as $L^\alpha$ regardless of the reset rate and the asymmetry of the initial condition, a behaviour completely different from that for the known case of Brownian Motion ($\alpha=2$).
For the symmetric $\alpha$-stable case with $\alpha \in (0, 2)$, we also found that taking $r \to \infty$ with fixed $L$ yields a finite MFPT, which therefore may be used as a valid approximation.
For completely asymmetric stable processes, we show that the scaling is very different from that for LF, and point out a method to obtain the exact distribution of FPT under resetting.
The analytical calculations are verified by numerical simulations of the processes under resetting.
\end{abstract}


\maketitle

%
%
\section{Introduction\label{sec:introduction}}

Brownian Motion (BM)~\cite{borodin2002handbook} is the simplest and very well studied model for diffusion and random motion. 
However, to describe many experiments, one needs more sophisticated models that capture some aspects of anomalous diffusion.
A common class of such stochastic processes are L\'evy flights (LF)~\cite{dubkov2008,chechkin2008introduction}, which can be used to describe various phenomena ranging from finance~\cite{laherrere1998,mantegna2000,lera2018gross} to the study of animal movement~\cite{brockmann2006,sims2008}.
Other significant classes of processes include L\'evy processes with only downward jumps, known as spectrally negative L\'evy processes (SNLPs). These processes are widely used in ruin theory, financial mathematics, and queuing theory.

A question commonly related to models of random motion is that of the first exit time from a given set or the first entrance to a given set.
These random times are often called \emph{first passage times} (FPTs)~\cite{redner2001, metzler2014firstpassage, grebenkov2024target}.
In the biological context it may be a protein exiting a chanel and in the financial setup we may ask when a price will exceed a given threshold. 
From the modeling point of view, due to the stochastic nature of the model, the FPT would be a random variable, and the most interesting features are its distribution and mean value (MFPT).

Stochastic (total) resetting~\cite{evans2011diffusion,evans2020stochastic,gupta2022stochastic} is a simple protocol that affects the original process by allowing it to restart.
At random times, called resetting events, the process is terminated and starts anew, losing its previous information.
The times between consecutive resetting events are usually independent of the process and drawn from a given probability distribution, the most popular being an exponential distribution~\cite{evans2011diffusion} (we call such resetting Poissonian).
Models that include stochastic resetting have various applications, ranging from search strategies~\cite{reynolds2009, viswanathan2011physics, palyulin2014levy} through chemistry~\cite{reuveni2014role} to computer science~\cite{luby1993optimal}.

The reason why stochastic resetting is important in the context of FPT problems is that it not only changes the process itself~\cite{bello2025partial, colantoni2025master, barkai2023ergodic}, but the FPT under resetting presents universal properties~\cite{reuveni2016optimal, pal2017first, chechkin2018random} and the mean first passage time (MFPT) can be altered by resetting~\cite{reuveni2016optimal, evans2011diffusion}.
In some models, the MFPT without resetting can be divergent, but stochastic resetting makes its value finite~\cite{evans2011diffusion}.
Moreover, if the distribution of times between consecutive resetting events can be tuned, stochastic resetting can be used as a tool for optimisation~\cite{radice2024optimizing}. 

In this letter, we will focus on the first exit time from the interval $(-L, L)$ of a L\'evy flight under total resetting.
Our main result is an analytical expression for the asymptotic behaviour of MFPT as $L \to \infty$.

In the case without resetting, the full distribution of FPT is known~\cite{borodin2002handbook} only for the Brownian case ($\alpha = 2$).
For the anomalous case $\alpha \in (0, 2)$ the mean and variance are known~\cite{blumenthal1961} since the 1960s, and the full distribution still remains an open problem.

For the model under stochastic resetting, one may use the general framework~\cite{pal2017first, chechkin2018random} to derive MFPT based on the distribution of FPT without resetting.
Therefore, the analytical formula for MFPT is available for the Gaussian case ($\alpha = 2$)~\cite{pal2019first}, but analogous results for LF with $\alpha \in (0, 2)$ cannot be obtained using these methods.

In the case $\alpha \in (0, 2)$ only partial results are available.
In~\cite{zbik2024levy} the authors analytically calculated the coefficient of variation~\cite{reuveni2016optimal, pal2017first}, showing that the initial positions of the particle $x_0$ for which there exists a resetting optimising the MFPT depends on the stability index $\alpha$.
Other known results usually refer to different models such as CTRW~\cite{colantoni2025master}, L\'evy walks~\cite{zbik2024levy} or other models~\cite{radice2024optimizing, evans2020stochastic}. 

The main results of this letter are as follows.
Firstly, under resetting, the asymptotic scaling of the MFPT for $\alpha < 2$ is qualitatively different from that for $\alpha = 2$.
This behaviour is not observed without resetting, where the MFPT asymptotically behaves like $\propto L^\alpha$ for all $\alpha \in (0,2]$ and large $L$.
Secondly, the scaling for $\alpha < 2$ is independent of the reset rate $r$ while for Brownian motion it depends on $r$.
Thirdly, numerical analysis strongly corroborates the theoretical results.  Finally, surprisingly, the case $r \to \infty$ yields a finite MFPT, which makes it a potentially very useful approximation.
We would like to emphasise that we also analyse SNLP.

The letter is organised as follows. The model and the quantities of interest are defined in Section~\ref{sec:model}.
We then prove the main result of this paper in Section~\ref{sec:results} and briefly discuss the findings.
The letter ends with a summary in Section~\ref{sec:summary}. 

%
%
\section{Model\label{sec:model}}

We consider a random walker that can be described as a L\'evy flight in 1 dimension or as a spectrally negative L\'evy process.
The position of such a particle is given by a Langevin equation
\[ \begin{cases}
\frac{\d Y(t)}{\d t} = \xi_\alpha \\
Y(0) = x_0,
\end{cases} \]
where $\xi_\alpha$ is either the symmetric $\alpha$-stable Lévy-type noise with $\alpha \in (0, 2]$ or, in the case of SNLP, appropriate L\'evy-type noise.
$Y(t)$ represents the position of the particle at time $t \in [0, \infty)$, and $x_0$ is its initial position.
The increments of such process $\Delta Y = Y(t + \Delta t) - Y(t)$ are independent, stationary, and in the case of LF distributed according to an $\alpha$-stable law whose characteristic function reads
\[ \avg{ \e^{i k \Delta Y}} = \exp ( - \Delta t | \sigma k |^\alpha ), \]
where $\avg{\cdot}$ denotes the expectation and $\sigma > 0$ is a scale parameter.
The case $\alpha = 2$ is the standard Brownian Motion, the case $\alpha \in (0, 2)$ corresponds to anomalous diffusion.

We then apply stochastic resetting~\cite{evans2011diffusion} to the process, which at random times brings the particle to its initial position $x_0$.
We will consider Poissonian resetting with reset rate $r$, meaning that the random times between consecutive resetting events are i.i.d.\,and distributed according to the exponential distribution whose density reads
\[ \phi(t) = \begin{cases}
r \e^{-r t}  \quad &\text{for}~t > 0 \\
0 \quad &\text{for}~t \le 0.
\end{cases} \]
The position of the particle under resetting will be denoted by $X(t)$.

Our quantity of interest are the first passage times (FPTs) of $Y$ and $X$, being the first (random) times when the particle 
evolving according to $Y$ or to $X$ escapes the interval $(-L, L)$
\begin{equation}
\label{eq:fpt-def}
\tau_L = \tau^{x_0}_L = \inf \left \{ t \ge 0 \colon Y(t) \notin (-L,~L) \right \}, 
\end{equation}
\begin{equation}
\T_L = \T_L^{x_0}=\inf \left \{ t \ge 0 \colon X(t) \notin (-L,~L) \right \}. 
\end{equation}
If it does not lead to disambiguation, we will skip the superscript $^{x_0}$ denoting the initial position of the process.

Then we may consider the mean first passage time (MFPT) defined as the average over trajectories 
\begin{equation*}
\avg{\tau_L} \quad \text{and} \quad \avg{\T_L}.
\end{equation*}
In this letter, we will analytically derive the asymptotic behaviour of $\avg{\T_L}$ when $L \to \infty$.

In the following derivations, by $\varphi_Z(s)$ we will denote the Laplace transform of a r.v.\,$Z \ge 0$ with the parameter $s$, i.e.
\[ \varphi_Z(s) \equiv \avg{\e^{-s Z}}. \]

%
%
\section{Results\label{sec:results}}
The general formula for MFPT under resetting reads
\begin{equation}\label{general}
\avg{ \T_L } = \frac{1}{r} \frac{1 - \varphi_{\tau_L}(r)}{\varphi_{\tau_L}(r)} , \end{equation}
see, for example,~\cite{reuveni2016optimal}.
In what follows, we will separately study LF and SNLP.

%
\subsection{L\'evy flights}

We first focus on the case where $Y$ is a LF.  
Due to the scaling properties of LF (law of $Y(At) - x_0$ is the same as of $A^\alpha (Y(t) - x_0)$ for $A \in \R_{> 0}$) the laws of $\tau_L^{x_0}$ and $L^\alpha \tau_1^{x_0 / L}=L^\alpha \tau_1$ are the same.
This allows us to simplify the formula for MFPT
\begin{equation}\label{formulamean}
\avg{ \T_L } = \frac{1}{r} \frac{1 - \varphi_{\tau_1}(r L^\alpha)}{\varphi_{\tau_1}(r L^\alpha)},
\end{equation}
and to find the asymptotic behaviour for large $L$, one only needs to know the behaviour of $\varphi_{\tau_1}(s)$ for large $s$.

However, for $\alpha \in (0, 2)$ we know the small time behaviour of FPT~\cite[Thm. 2.1]{10.3150/13-BEJ517}
\begin{equation}\label{smalltime}
\p( \tau_1 \le \varepsilon) \sim A \varepsilon \quad \text{as}\;\varepsilon \to 0,
\end{equation}
where 
\[A = \int_{1 - \frac{x_0}{L}}^\infty g(x)\d x +  \int_{1 + \frac{x_0}{L}}^\infty g(x)\d x,\]
and 
\[g(x) = 
\begin{cases}
\begin{aligned}
\frac{\alpha \sigma^\alpha}{2 \Gamma (1-\alpha) \cos \left( \frac{\pi \alpha}{2} \right) } &\cdot \frac{1}{|x|^{1 + \alpha}}
&\quad \text{for}~\alpha \ne 1 \\
\frac{\sigma}{\pi} &\cdot \frac{1}{x^2}
&\quad \text{for}~\alpha = 1
\end{aligned}
\end{cases} \]
is the density of the L\'evy measure of the LF $Y$.
Evaluating the integrals for $\alpha \ne 1$ yields 
\[ A = \frac{\sigma^\alpha }{2 \Gamma (1-\alpha ) \cos \left(\frac{\pi  \alpha }{2}\right)} \cdot 
\left( (1- \frac{x_0}{L})^{-\alpha} + (1 + \frac{x_0}{L})^{-\alpha} \right). \]
If the initial condition is symmetric ($x_0 = 0$) it simplifies to 
\[ A = \frac{\sigma^\alpha }{\Gamma (1-\alpha ) \cos \left(\frac{\pi  \alpha }{2}\right)}. \]
Formally $\alpha = 1$ requires a separate calculation, but taking the limit as $\alpha \to 1$ gives the same result
\[ \frac{\sigma}{\pi} \cdot \left( (1- \frac{x_0}{L})^{-\alpha} + (1 + \frac{x_0}{L})^{-\alpha } \right). \]

\begin{figure}[h]
\includegraphics[width=\columnwidth]{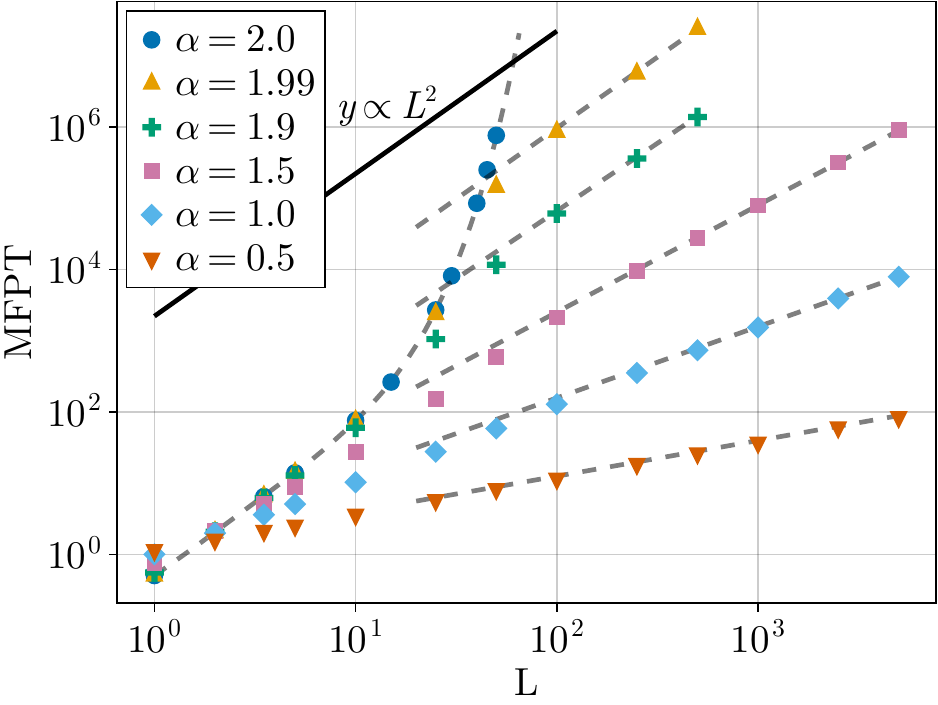}
\caption{MFPT as a function of the interval half-length $L$ plotted in a log-log scale. For $\alpha < 2$ the dashed lines are power-laws given by~\eqref{eq:final-scaling}, and for $\alpha = 2$ the known analytical expression for MFPT. The solid line visualises the inclination of a hypothetical $L^2$ power law. Parameters used for the simulation: $r = 0.05$, $\Lambda = 0.0$ and $\sigma = 1.0$.
}
\label{fig:power-laws}
\end{figure}

Then by a Tauberian theorem~\cite[XIII.5 Theorem 2]{feller1971vol2}, one may relate the asymptotic of $\varphi_{\tau_1}$ at infinity with the distribution of $\tau_1$ close to~$0$.
Note that this procedure is not necessary for the case of BM ($\alpha = 2$), where the exact formula for MFPT is known~\cite[Eq. (15)]{pal2019first}
\[
\avg{ \mathcal{T}_L } = \frac{1}{r} \left[ \frac{\sinh \frac{2L}{\sqrt{\sigma^2/r}} }{ \sinh \frac{L-x_0}{\sqrt{\sigma^2/r}} + \sinh \frac{x_0+L}{\sqrt{\sigma^2/r}} } - 1 \right].
\]
In the above formula $\sigma = 1/\sqrt{2}$ corresponds to the Wiener process.
Applying a Tauberian theorem (in case BM just using the exact formula) gives, if the initial condition is symmetric, i.e.\ $x_0 = 0$,  the final scaling of MFPT for $L \to \infty$
\begin{equation}\label{eq:final-scaling}
\avg{\T_L} \sim 
\begin{cases}\begin{aligned}
\cos \left(\frac{\pi  \alpha }{2}\right) \Gamma (1-\alpha ) ~ \frac{L^\alpha}{\sigma^\alpha} \quad &\text{for}~\alpha \in (0, 2) \\
\frac{1}{r} ~ \exp{( \sqrt{r / \sigma^2} \cdot L )} \quad &\text{for}~\alpha = 2,
\end{aligned}\end{cases} 
\end{equation}
where $\Gamma(\cdot)$ is the Gamma function. 
If the initial condition is not symmetric $(x_0 \ne 0)$, then the scaling also depends on the ratio $\Lambda := \frac{x_0}{L}$ and reads
\begin{equation}\label{eq:final-scaling-asymetric}
\avg{ \T_L } \sim 
\begin{cases}\begin{aligned}
C\frac{L^\alpha}{\sigma^\alpha} \quad &\text{for}~\alpha \in (0, 2) \\
\frac{1}{r} ~ \exp{(\gamma \cdot L)} \quad &\text{for}~\alpha = 2,
\end{aligned}\end{cases}
\end{equation}
where 
\begin{align*}
C &= \frac{2 \cos \left(\frac{\pi \alpha}{2}\right) \Gamma (1-\alpha )}{
 (1 - \Lambda)^{-\alpha }+ (1 + \Lambda )^{-\alpha } },\\
\gamma &= \sqrt{r / \sigma^2} \min \{1 - \Lambda, 1 + \Lambda\}.
\end{align*}
For $\alpha = 1$ the formula needs to be treated as a limit $\alpha \to 1$ to avoid $0 \cdot \infty$
\[ \avg{ \T_L } \sim \frac{\pi}{(1 - \Lambda)^{-1 }+ (1 + \Lambda )^{-1 } }~ \frac{L}{\sigma} . \]
The analytical derivations presented above agree with Monte Carlo simulations, as shown in FIG~\ref{fig:power-laws}.

We would like to stress that the asymptotic derived in equation~\eqref{eq:final-scaling} is a rare example of an analytical result for LF where most of the findings are numerical.
Note the qualitative change in the scaling of MFPT when we go from the gaussian to $\alpha$-stable case (an exponent changes to a power law), which is not observed without resetting~\cite{blumenthal1961, dybiec2017levy}, since then MFPT reads
\begin{equation}\label{eq:mfpt-no-reset}
\avg{ \tau_L} = \frac{(1-\Lambda^2)^{\alpha/2}}{\Gamma(1+\alpha) } \times  \frac{L^\alpha}{\sigma^\alpha} 
 \quad \text{for} ~ \alpha \in (0, 2]. 
\end{equation}
This qualitative change comes from the fundamentally different nature of LF and BM.
For $\alpha < 2$, due to long-range jumps, the trajectories are discontinuous, which is not the case for BM.
When $L$ is large, a Brownian particle is almost prevented from exploring far, as resetting keeps it near the origin.
However, LF, even under resetting, may still experience long jumps, which give the main contribution to the FPT.

\begin{figure}[h]
\includegraphics[width=\columnwidth]{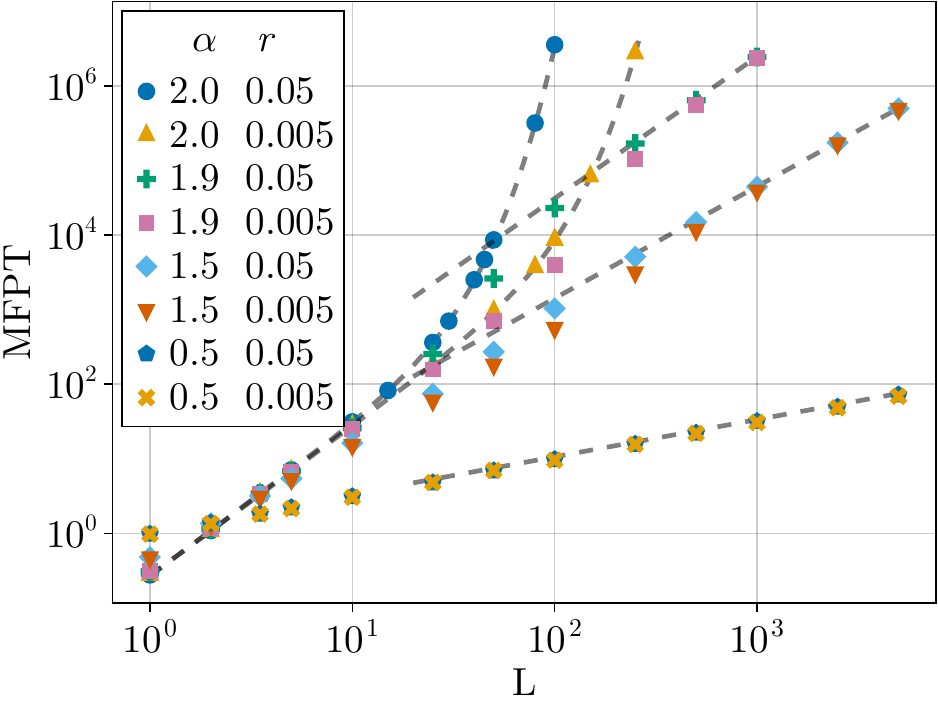}
\caption{MFPT as a function of the interval half-length $L$ plotted in a log-log scale. For $\alpha < 2$ the dashed lines are power-laws given by~\eqref{eq:final-scaling-asymetric}, and for $\alpha = 2$ the analytical expression for MFPT.
The parameters used were $\Lambda = 0.3$, and $\sigma = 1.3$.
}\label{fig:different-r}
\end{figure}

Interestingly, the power-law scaling for $\alpha < 2$ is independent of the reset rate $r$, which is not the case for BM; see FIG.~\ref{fig:different-r}.
This phenomenon is once again related to the presence of long jumps for LF.
As the escape is usually related to a single big jump, the resetting events do not hinder the MFPT. 

\begin{figure}[h]
\includegraphics[width=\columnwidth]{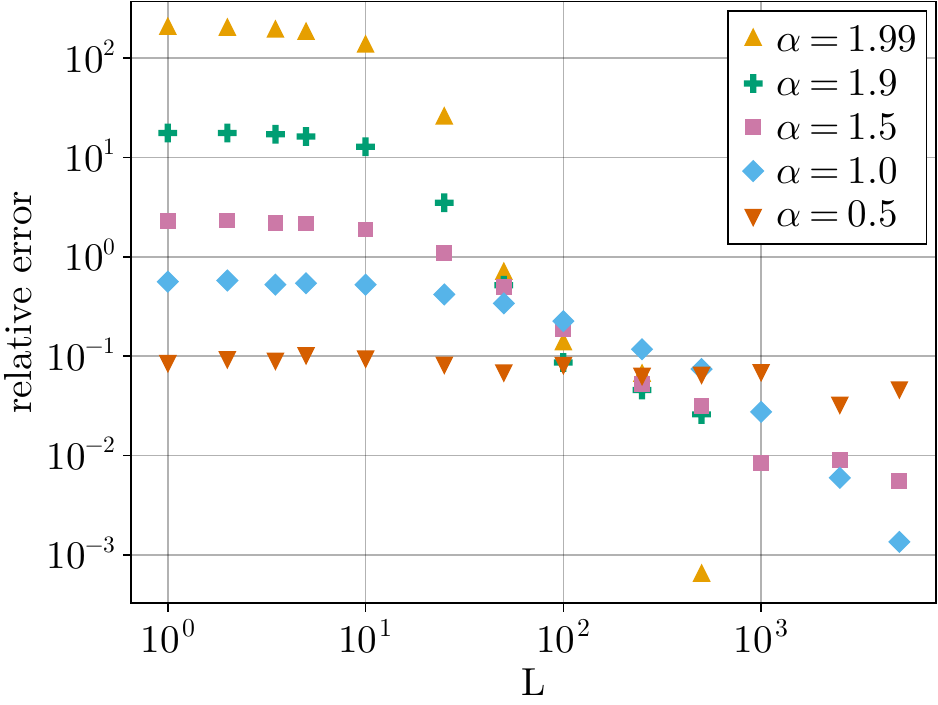}
\caption{Relative error between the large $L$ asymptotic of MFPT, formulas~\eqref{eq:final-scaling}-\eqref{eq:final-scaling-asymetric}, and Monte Carlo simulations.
The relative error was calculated as the absolute value of the difference of MFPT coming from numerical simulations and the asymptotic formula, all divided by the value of MFPT from the simulations. The parameters used were $\Lambda = 0.0, \sigma = 1.0$, and $r = 0.05$.
}\label{fig:errors}
\end{figure}

Based on the numerical simulations, we analysed the quality of approximating MFPT through formulas~\eqref{eq:final-scaling} and~\eqref{eq:final-scaling-asymetric}, see Fig.~\ref{fig:errors}.
For the sets of tested parameters, $r \avg{ \T_1} \approx 10^{-2} - 10^{-1}$, the error dropped below 10\% at the scale $L/\sigma \approx 10^2 - 10^3$ and decreased further with increasing $L$.
Figure~\ref{fig:errors} is a good representative of the observed behaviours.
For small $L$ the approximation, as expected, does not work and improves as $L$ grows.
The error decreases more rapidly for higher values of $\alpha$.
We used the Euler--Maruyama discretisation technique~\cite{higham2001algorithmic, mannella2002} to simulate the Langevin equation
\[ x_{n+1} = x_n + (\Delta t)^{1/ \alpha} \xi_{n+1} , \]
where $\xi_i$ are i.i.d. random variables drawn from the symmetric $\alpha$-stable distribution.
The exact time steps used for the simulation differed according to $L$ and $r$ to match the temporal and special scales of the problem.
For the data in Figures~\ref{fig:power-laws}--\ref{fig:errors} the ratio of MFPT and $\Delta t$ was at least of order $10^2$ and for $\alpha \in (0.75, 2.0)$ above $10^3$.
Each point corresponds to a mean taken from $10^3$ or $10^4$ trajectories.

\begin{figure}[h]
\includegraphics[width=\columnwidth]{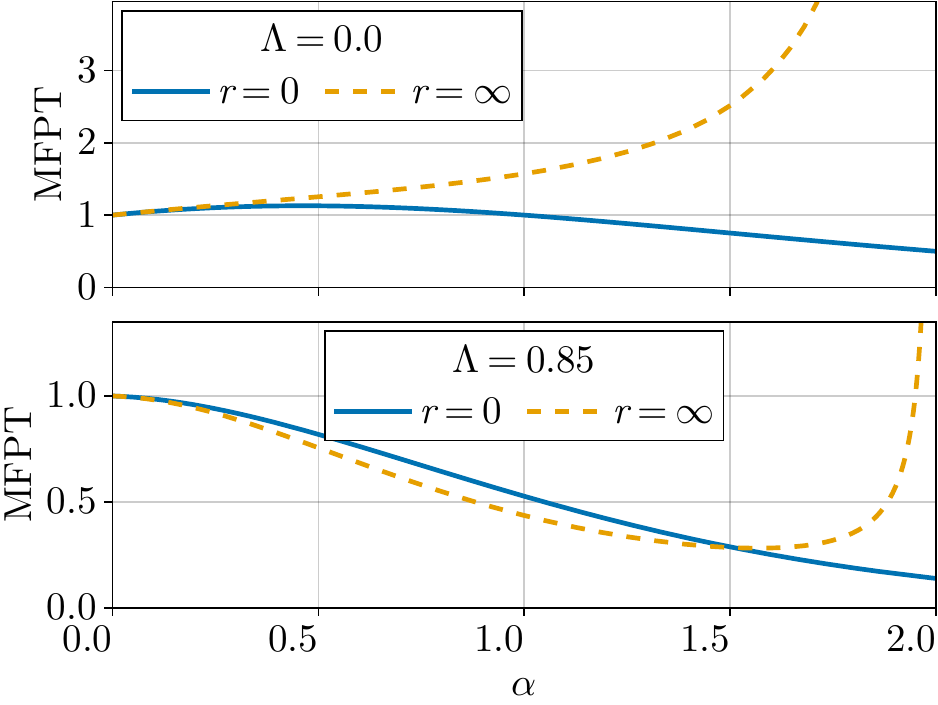}
\caption{MFPT for LF escaping an interval without resetting (solid, blue) or with the reset rate $r \to \infty$ (dashed, orange).
In the top panel the initial condition is symmetric $\Lambda = 0.0$, while the bottom panel corresponds to $\Lambda = 0.85$.
In both cases $L = 1$ and $\sigma = 1$.
}\label{fig:r-infty}
\end{figure}

We also note that instead of taking the limit $L \to \infty$, one may consider the case where $L$ is fixed and $r \to \infty$.
The obtained asymptotic formulas turn out to be the same as previously, Eqs.~\eqref{eq:final-scaling} and~\eqref{eq:final-scaling-asymetric}, due to representation~\eqref{formulamean}.
Therefore, one may directly compare Eq.~\eqref{eq:final-scaling-asymetric} with the formula for MFPT without resetting~\eqref{eq:mfpt-no-reset}, see Figure~\ref{fig:r-infty}.
This corresponds to studying a LF that escapes an interval (with fixed $L, \sigma, x_0$), where we compare the case without resetting to the case where resetting appears infinitely often.
Firstly, it turns out that in this scenario, the MFPT under resetting for $\alpha < 2$ does not explode.
In addition, the dependence on the asymmetry of the initial condition $\Lambda$ is different based on the presence of resetting. 
Surprisingly, there is a parameter regime, e.g. $\Lambda = 0.85$ and $\alpha = 1$, for which $\avg{ \T_L} < \avg{ \tau_L }$.
That is, even resetting with very large intensity may still lower the MFPT compared to the case without resetting.

%
\subsection{Spectrally negative L\'evy processes}

We can also consider the case where the process $Y$ without resetting does not have upward jumps, that is, $Y$ is a SNLP.
In this case, by~\cite[Thm. 8.8.1]{kyprianou2014fluctuations} 
\[ \varphi_{\tau_L}(r) = (1-Z^{(r)}(2L))\frac{W^{(r)}(x_0 + L)}{W^{(r)}(2L)}+Z^{(r)}(x_0 + L),\]
where $W^{(r)}(x)$ is the so-called scale function with the Laplace transform $1/(\varphi_{Y(1)}(s)-r)$ and
$Z^{(r)}(x)=1+r\int_0^x W^{(r)}(y)\d y$.

Let us denote the inverse of $\varphi_{Y(1)}$ by $\Phi(r)$, that is $\varphi_{Y(1)}(\Phi(r))=r$.
Then the asymptotic behaviour of $W^{(r)}(x)$ and $Z^{(r)}(x)$ as $x\to+\infty$ read (see~\cite[p. 251-252]{kyprianou2014fluctuations})
\[ W^{(r)}(x) \sim De^{\Phi(r) x} \quad \text{and} \quad  Z^{(r)}(x) \sim \frac{r}{\Phi(r)}De^{\Phi(r) x} \]
where $D= [\varphi_{Y(1)}'(\Phi(r))]^{-1}$.
Hence in the case of SNLP, by~\eqref{general}, we get
\[\avg{ \T_L }\sim \frac{1}{r}\exp(\Phi(r)(L-x_0)).\]

Finally, if $Y$ is a completely asymmetric stable process with $\alpha \in (1, 2]$, then by~\cite[Thm. 1]{Bertoinexit}
\begin{align*}
\varphi_{\tau_L}(r) &= E_\alpha(r[L+x_0]^\alpha) \\
&\quad - \frac{[L + x_0]^{\alpha - 1}}{2^{\alpha - 1} L^{\alpha - 1}}
\frac{E_\alpha'(r[L+ x_0]^\alpha)}{E_\alpha'(r2^\alpha L^\alpha)}(E_\alpha(r2^\alpha L^\alpha)-1),
\end{align*}
where $E_\alpha=\sum_{n=0}^\infty x^n/\Gamma(1+\alpha n)$ is a Mittag-Leffler function of parameter $\alpha$.
Note, that for $\alpha = 2$ we recover the BM case, as $E_2(x)=\cosh(\sqrt{x})$ and $E_2'(x)=\frac{\sinh (\sqrt{x})}{2\sqrt{x}}$.
Now knowing the asymptotic behaviour of $E_\alpha(x)\sim \frac{1}{\alpha}\exp(x^{1/\alpha})$ and $E_\alpha'(x)\sim \frac{1}{\alpha^2}x^{(1-\alpha)/\alpha}\exp(x^{1/\alpha})$ as $x \to +\infty$, we find, via~\eqref{general}, that for $\Lambda \ge 0$
\begin{equation*}
\avg{\T_L} \sim \frac{1}{r}\exp(r^{1/\alpha} [1 - \Lambda] L).
\end{equation*}
What is worth noting is that for $\alpha \in (1, 2)$, the asymptotic behaviour for completely asymmetric stable processes is very different from that for the symmetric case.
For $\alpha = 2$, we also recover the result for BM~\eqref{eq:final-scaling}.
Finally, it is important that via~\eqref{general}, one can also get the exact expression for MFPT under resetting expressed in terms of Mittag-Leffler functions and their derivatives.
Similarly, using~\cite[Eq.~(6)]{reuveni2016optimal} instead of~\eqref{general}, the full Laplace transform of $\T_L$ can be obtained.

%
%
\section{Summary and conclusions\label{sec:summary}}

In this letter, we studied the problem of first passage time for L\'evy flights and spectrally negative L\'evy processes under stochastic resetting escaping an interval $(-L, L)$.
When the stability index $\alpha = 2$ (Brownian Motion), the full distribution of FPT is known both with and without resetting.
However, for LF with $\alpha \in (0, 2)$ and without resetting, only the mean and variance of FPT are known, with no such analytical results for the case under resetting.

The main result of this letter is the asymptotic formula, Eqs.~\eqref{eq:final-scaling}-\eqref{eq:final-scaling-asymetric}, for MFPT for LF under stochastic resetting when $L \to \infty$.
We found that introducing resetting qualitatively changes the scaling for the BM case, a power law $L^2$ changes to an exponent $\e^{A L}$, whereas for $\alpha \in (0, 2)$ we observe a power law behaviour $L^\alpha$ regardless of the resetting.
Moreover, for $\alpha \in (0, 2)$, the asymptotic value of MFPT is independent of the reset rate $r$ and only the proportionality factor depends on the asymmetry of the initial condition $\Lambda = \frac{x_0}{L}$.
This differs significantly from the behaviour known for BM, where the coefficient in the exponent depends on both $r$ and $\Lambda$, making BM much more responsive to the presence of resetting.

If instead of $L \to \infty$ one considers $r \to \infty$ with fixed $L$, the derived formulas, Eqs.~\eqref{eq:final-scaling}-\eqref{eq:final-scaling-asymetric}, are still correct.
This results in an interesting setup, for which despite infinitely often resetting, we still observe finite MFPT.
Moreover, by applying such resetting, MFPT can even become smaller than without resetting.
Therefore, this a priori unphysical case may still be applicable with $r \to \infty$ treated as an approximation.

Finally, we also studied the case when the L\'evy process does not have positive jumps (spectrally negative L\'evy process).
Analogously to the case of LF, we derived the asymptotic for MFPT under resetting as $L \to +\infty$.
For the case of completely asymmetric stable processes with $\alpha \in (1, 2]$, we explicitly calculated the scaling parameters, reviling a completely different scaling behaviour than for LF.
For completely asymmetric stable processes, we also pointed out how to find the full distribution of FPT and the exact formula for the MFPT under resetting in terms of Mittag-Leffler functions.

%
%
\section*{Acknowledgments}
The research is partially supported by the Polish National Science Centre
Grant No. 2023/51/B/ST1/01270.

The authors thank Jagiellonian University for its support and resources during the early stages of this work.

\section*{Data availability}
Data (generated randomly using the model presented in the article) supporting the findings of this study are available from the corresponding author (B\.Z) upon request. 

%
%

\section*{References}
\bibliography{core-bibliography}

\end{document}